\documentclass[11pt]{article}
\usepackage{amssymb}
\usepackage{amsthm}
\usepackage{graphicx}
\usepackage{graphics}
\usepackage{amsthm}

\input pictex.tex
\newcommand{\esp}{\hspace{0.01cm}}

\theoremstyle{definition}

\newcommand{\ce}{\mathrm{C}}

\input epsf

\begin{document}

\date{}
\author{Andr\'es Navas}
%University of Santiago}

\title{An example concerning the Theory of Levels for codimension-one foliations}
\maketitle

%\noindent{\bf Abstract.} 

%\vspace{-1cm}

\noindent{\Large {\bf Introduction}}

\vspace{0.5cm}

An important aspect of foliations concerns the existence of local minimal sets. Recall  
that a foliated manifold has the LMS property if, for every open, saturated set $W$ and 
every leaf $L \subset W$, the relative closure $\bar{L} \cap W$ contains a minimal set 
of $F|_W$. A fundamental result (due to Cantwell-Conlon \cite{CC} and Duminy-Hector 
\cite{hector}) establishes that for codimension-one foliations which are transversely 
of class $C^{1 + \mathrm{Lipschitz}}$, the LMS property holds. This is the basic 
tool of the so-called {\em Theory of Levels}. 

A well known example due to Hector (which corresponds to the suspension of a group action 
on the interval) shows that the LMS property is no longer true for codimension-one 
foliations which transversely are only continuous (see \cite{CC-libro}, Example 8.1.13). 
Despite this, in the recent years the possibility of extending some of the aspects of 
the Theory of Levels to smoothness smaller than $C^{1 + \mathrm{Lipschitz}}$ has 
naturally appeared $\cite{CC2,DKN}$. In this Note we will show that, however, 
analogues of Hector's example appear in class $C^1$ (and actually in class 
$C^{1+\alpha}$ for some values of $\alpha$).

%%%%%%%%%%%%%%%%%%%%%%%%%%%%%%%%%%%%%%%%%%%%%%%%%%%%%%%%%%%%%%%%%%%%%%%%%%%%%%%%%%%%%

\section{A General Construction}

\hspace{0.45cm} Let $(a_n)_{n \in \mathbb{Z}}$ 
be a sequence such that $a_{n+1} < a_{n}$ for all 
$n \in \mathbb{Z}$, $a_n \to 0$ as $n \to \infty$, and $a_n \to 1$ as $n \to -\infty$. 
Let $(n_k)$ be a strictly increasing sequence of positive integers, and let 
$f \!: [0,1] \rightarrow [0,1]$ be a homeomorphism such that $f(a_{n+1}) = a_n$ 
for all $n \in \mathbb{Z}$. For each $k > 0$ we let $u_k,v_k,b_k,c_k$ be such that 
$a_{n_k+1}<b_k<u_k<v_k<c_k<a_{n_k}$. For each $i \in \{0,\ldots,n_{k+1}-n_k\}$ 
we define $u_k^i = f^i (u_k)$ and $v_k^i = f^i (v_k)$. Notice that
$$f^i([u_{k+1}^0,v_{k+1}^0]) = [u_{k+1}^i,v_{k+1}^i] \subset 
f^i ([a_{1+n_{k+1}},a_{n_{k+1}}]) = [a_{n_{k+1}-i+1},a_{n_{k+1}-i}].$$ 

\vspace{0.1cm}

Now we let $g \! : [0,1] \rightarrow [0,1]$ be a homeomorphism such that:\\

\noindent -- $g = Id$ on $[a_{n+1},a_n]$ for each $n < 0$ 
and each $n > 0$ such that $n \neq n_k$ for every $k$,

\noindent -- $g = Id$ on $[a_{1+n_k},b_k] \cup [c_k,a_{n_k}]$, \esp 
$g(u_k^0) = v_k^0$, \esp and $g$ has no fixed point on $]b_k,c_k[$.

\vspace{0.25cm}

\noindent{\bf Main assumption:} In order that $f,g$ generate a group of 
homeomorphisms of $[0,1]$ whose associated suspension does not have  
the LMS property, we assume that (see Figure 1)
$$u_{k+1}^{n_{k+1} - n_k} = b_k \qquad \mbox{ and } 
\qquad v_{k+1}^{n_{k+1} - n_k} = c_k.$$

With these general notations, Hector's example  corresponds to the choice $n_k = k$. We 
will show that, by taking $n_k = 2^k$, one may perform this construction in such a way 
the resulting maps $f$ and $g$ are diffeomorphisms of class $\ce^1$ (actually, 
of class $\ce^{1 + \alpha}$ for any $\alpha < (\sqrt{5} - 1) / 2$). It is quite 
possible that improving slightly the method one can smooth the action up to the class 
$\ce^{2-\delta}$ for any $\delta > 0$; compare \cite{tsuboi}, where for a similar  
construction T. Tsuboi deals with the $\ce^{3/2 - \delta}$ case before the 
$\ce^{2-\delta}$ case due to technical difficulties.

\vspace{2.4cm}

%%%%%%%%%%%%%%%%%%%%%%%%%%%%%%%%%%%%%%%%%%%%%%%%%%%%%%%%%%%%%%%%%%%%%%%%%%%%%%%%%%%%%%%%%%%%%%
\beginpicture

\setcoordinatesystem units <0.9cm,0.9cm>

\putrule from 0 0 to 6 0
\putrule from 10 0 to 16 0

\begin{small}
\putrule from 2.25 0.01  to 3.87 0.01  
\putrule from 2.25 -0.01 to 3.87 -0.01
\putrule from 2.25 0.02  to 3.87 0.02  
\putrule from 2.25 -0.02 to 3.87 -0.02
\put{$b_{k+1}$} at 1.2 -0.6 
\put{$c_{k+1}$} at 4.8 -0.6 
\put{$($} at 1.2 0 
\put{$)$} at 4.8 0 

\putrule from 11.18 0.01  to 14.825 0.01  
\putrule from 11.18 -0.01 to 14.825 -0.01
\putrule from 11.18 0.02  to 14.825 0.02  
\putrule from 11.18 -0.02 to 14.825 -0.02
\put{$b_k$} at 11.2 -0.6 
\put{$c_k$} at 14.8 -0.6 
\put{$($} at 11.2 0 
\put{$)$} at 14.8 0 
\end{small}

\begin{Large}
\put{$a_{1+n_{k+1}}$} at -0.4 -0.9
\put{$a_{n_{k+1}}$} at 5.6 -0.9 
\put{$($} at -0.4 0 
\put{$)$} at 5.6 0 

\put{$f^{n_{k+1} - n_{k}}$} at 7.65 1.75 
\put{$g$} at 12.7 1 
\put{$g$} at 2.7 1 

\put{$a_{1+n_{k}}$} at 9.6 -0.9
\put{$a_{n_{k}}$} at 15.6 -0.9 
\put{$($} at 9.6 0 
\put{$)$} at 15.6 0 
\end{Large}

\begin{tiny}
\put{$u_{k+1}$} at 1.7 -0.35
\put{$v_{k+1}$} at 3.3 -0.35 
\put{$($} at 1.7 0
\put{$)$} at 3.3 0 

\put{$u_k$} at 11.7 -0.35
\put{$v_k$} at 13.3 -0.35 
\put{$($} at 11.7 0
\put{$)$} at 13.3 0 
\end{tiny}

\circulararc 45 degrees from 11.25 1.5  
center at 7.25 -8

\plot 
11.25 1.5 
11.03 1.5 / 

\plot 
11.25 1.5 
11.085 1.67 /

\circulararc 180 degrees from 13.12 0.7  
center at 12.35 0.7 

\plot 
13.12 0.7 
13.18 0.92 / 

\plot 
13.12 0.7 
13 0.9 /

\circulararc 180 degrees from 3.12 0.7  
center at 2.35 0.7 

\plot 
3.12 0.7 
3.18 0.92 / 

\plot 
3.12 0.7 
3 0.9 /

\put{Figure 1} at 7.12 -1.5

\put{} at -1.4 0

\endpicture

%%%%%%%%%%%%%%%%%%%%%%%%%%%%%%%%%%%%%%%%%%%%%%%%%%%%%%%%%%%%%%%%%%%%%%%%%%%%%%%%%%%%%%%%%%%%%%

\vspace{0.75cm}

%%%%%%%%%%%%%%%%%%%%%%%%%%%%%%%%%%%%%%%%%%%%%%%%%%%%%%%%%%%%%%%%%%%%%%%%%%%%%%%%%%%%%%%%%%%%

\section{The length of the intervals and the estimates}

We let \esp $|[u_{k+1}^i,v_{k+1}^i]| = \lambda_k^i \esp |[u_{k+1},v_{k+1}]|$, 
\esp where the constant $\lambda_k \!>\! 1$ satisfies the compatibility relation 
\begin{equation}
\lambda_k^{2^k} = \frac{|[u_{k+1}^{2^k},v_{k+1}^{2^k}]|}{|[u_{k+1},v_{k+1}]|} 
= \frac{|[b_k,c_k]|}{|[u_{k+1},v_{k+1}]|}.
\label{compa}
\end{equation}
Let $\varepsilon > 0$ be very small (to be fixed in a while). 
We put:\\

\noindent -- $|[a_{n+1},a_n]| = \frac{c_{\varepsilon}}{(1 + |n|)^{1 + \varepsilon}}$, \esp 
where $c_{\varepsilon}$ is chosen so that $\sum_{n \in \mathbb{Z}} |[a_{n+1},a_n]| = 1$;\\

\noindent -- $|[b_k,c_k]| = \frac{1}{2} |[a_{2^k + 1}, a_{2^k}]| = 
\frac{c_{\varepsilon}}{2 (1 + 2^{k})^{1 + \varepsilon}}$, \esp where $k > 0$;\\

\noindent -- $|[u_k,v_k]| = |[b_k,c_k]|^{1 + \theta}$.

We also assume that the center of $[a_{2^k + 1}, a_{2^k}]$ coincides with the center of 
$[b_k,c_k]$ and with that of $[u_k,v_k]$. Moreover, for each $i \!\in\! \{0,\ldots,2^k\}$, 
the centers of $[u_{k+1}^i,v_{k+1}^i]$ and $[a_{2^{k+1} - i +1},a_{2^{k+1} - i}]$ do coincide. 
For the estimates concerning regularity we will strongly use the following lemma from \cite{growth}.

\vspace{0.53cm}

\noindent{\bf Lemma.} {\em Let $\omega: [0,\eta] \rightarrow [0,\omega(\eta)]$ 
be a modulus of continuity such that the function $s \mapsto s / \omega(s)$ 
is non increasing. If $I,J$ are closed non degenerate intervals such that 
\esp $1/2 \leq |I| / |J| \leq 2$ \esp and}
$$\left| \frac{|J|}{|I|} - 1 \right| \frac{1}{\omega(|I|)} \leq M,$$
{\em then there exists a $\ce^{1 + \omega}$ diffeomorphism $f\!: I \rightarrow J$
which is tangent to the identity at the end points and whose derivative has 
$\omega$-norm bounded from above by $6 \pi M$.} 

\vspace{0.53cm}

Actually, for $I \!=\! [a,b]$ and $J \!=\! [a',b']$ one may take 
$f = \varphi_{a',b'}^{-1} \circ \varphi_{a,b}$, where $\varphi_{a,b}$ 
is defined by (a similar definition stands for $\varphi_{a',b'}$)
$$\varphi_{a,b} (x) = 
-\frac{1}{(b-a)} \mathrm{ctg} \left( \pi \Big( \frac{x-a}{b-a} \Big) \right).$$
The condition on the derivative at the end points allows us to fit together the maps in 
order to create a diffeomorphism of a larger interval. Actually, if all of the involved 
sub-intervals of type $I,J$ satisfy the hypothesis of the Lemma above
for the same constant $M$, 
then the derivative of the induced diffeomorphism has $\omega$-norm bounded by $12 \pi M$. 

In what follows we will deal with the modulus of continuity $\omega(s) \!=\! s^{\alpha}$ 
for the derivative, where $\alpha > 0$. A constant depending on the three parameters 
$\alpha, \theta, \varepsilon$, and whose value is irrelevant for our purposes,  
will be generically denoted by $M$.

\vspace{0.25cm}

%%%%%%%%%%%%%%%%%%%%%%%%%%%%%%%%%%%%%%%%%%%%

\noindent{\bf Estimates for $f$:} The diffeomorphisms $f$ is constructed by fitting together 
maps sending (see Figure 2): 

\vspace{0.1cm}

\noindent (i) $[u_{k+1}^i,v_{k+1}^i]$ \esp\esp into 
\esp\esp $[u_{k+1}^{i+1},v_{k+1}^{i+1}]$,

\vspace{0.1cm}

\noindent (ii) $[a_{2^{k+1} - i}, u_{k+1}^i]$ \esp\esp into 
\esp\esp $[a_{2^{k+1} - i - 1},u_{k+1}^{i+1}]$,

\vspace{0.1cm}

\noindent (iii) $[v_{k+1}^i,a_{2^{k+1} - i - 1}]$ \esp\esp  into 
\esp\esp $[v_{k+1}^{i+1},a_{2^{k+1} - i - 2}]$.

\vspace{0.1cm}

For (i) we have
$$\left| \frac{|[u_{k+1}^{i+1},v_{k+1}^{i+1}]|}{|[u_{k+1}^{i},v_{k+1}^{i}]|} - 1 \right| 
\frac{1}{|[u_{k+1}^{i},v_{k+1}^{i}]|^{\alpha}} 
= |\lambda_k - 1| \frac{1}{(\lambda_k^{i} |[u_{k+1}^0,v_{k+1}^0]|)^{\alpha}} 
\leq |\lambda_k - 1| \frac{1}{|[b_{k+1},c_{k+1}]|^{(1 + \theta)\alpha}}.$$
Now from (\ref{compa}) one obtains
$$\lambda_{k}^{2^k} = \frac{\frac{c_{\varepsilon}}{2 (1 + 2^k)^{1+\varepsilon}}}
{(\frac{c_{\varepsilon}}{2 (1+2^{k+1})^{1 + \varepsilon}})^{1+\theta}} 
\leq M \Big( \frac{(1 + 2^{k+1})^{1 + \theta}}{1 + 2^k} 
\Big)^{1 + \varepsilon} \leq M 2^{k \theta (1 + \varepsilon)}.$$
From the inequality $|2^{\alpha} - 1| \leq \alpha$ 
(which holds for $\alpha$ positive and small) 
one concludes that
$$|\lambda_k - 1| \leq M \frac{k}{2^k}.$$
On the other hand,
$$\frac{1}{|[b_{k+1},c_{k+1}]|} \leq M (1 + 2^{k+1})^{1 + \varepsilon} 
\leq M 2^{k(1 + \varepsilon)}.$$
Therefore,
\begin{equation}
\left| \frac{|[u_{k+1}^{i+1},v_{k+1}^{i+1}]|}{|[u_{k+1}^{i},v_{k+1}^{i}]|} - 1 \right| 
\frac{1}{|[u_{k+1}^{i},v_{k+1}^{i}]|^{\alpha}} \leq 
M \frac{k}{2^k} 2^{k (1+\varepsilon) (1+\theta) \alpha}.
\label{primera}
\end{equation}

\vspace{3cm}

%%%%%%%%%%%%%%%%%%%%%%%%%%%%%%%%%%%%%%%%%%%%%%%%%%%%%%%%%%%%%%%%%%%%%%%%%%%%%%%%%%%%%%%%%%%%%%
\beginpicture

\setcoordinatesystem units <0.9cm,0.9cm>

\putrule from 0 0 to 6 0
\putrule from 10 0 to 16 0

\begin{small}
\putrule from 2.25 0.01  to 3.87 0.01  
\putrule from 2.25 -0.01 to 3.87 -0.01
\putrule from 2.25 0.02  to 3.87 0.02  
\putrule from 2.25 -0.02 to 3.87 -0.02
%\put{$b_{k+1}$} at 1.2 -0.6 
%\put{$c_{k+1}$} at 4.8 -0.6 
%\put{$($} at 1.2 0 
%\put{$)$} at 4.8 0 

\putrule from 12.25 0.01  to 13.87 0.01  
\putrule from 12.25 -0.01 to 13.87 -0.01
\putrule from 12.25 0.02  to 13.87 0.02  
\putrule from 12.25 -0.02 to 13.87 -0.02

\put{$A$} at 3.08  -0.35 
\put{$C$} at 13.08 -0.35 

\put{$B$} at 3.08  0.86 
\put{$D$} at 13.08 0.86 

\putrule from 0 0.86 to 2.85 0.86  
\putrule from 3.3 0.86 to 5.8 0.86 

\putrule from 10 0.86 to 12.85 0.86  
\putrule from 13.3 0.86 to 15.8 0.86   

\end{small}

\begin{Large}
\put{$a_{2^{k+1}-i}$} at -0.4 -0.6
\put{$a_{2^{k+1}-i-1}$} at 5.6 -0.6 
\put{$($} at -0.4 0 
\put{$)$} at 5.6 0 

\put{$f$} at 7.65 1.75 

\put{$a_{2^{k+1}-i-1}$} at 9.6 -0.6
\put{$a_{2^{k+1}-i-2}$} at 15.6 -0.6 
\put{$($} at 9.6 0 
\put{$)$} at 15.6 0 
\end{Large}

\begin{tiny}
\put{$u_{k+1}^i$} at 1.7 -0.35
\put{$v_{k+1}^i$} at 3.3 -0.35 
\put{$($} at 1.7 0
\put{$)$} at 3.3 0 

\put{$u_{k+1}^{i+1}$} at 11.7 -0.35
\put{$v_{k+1}^{i+1}$} at 13.3 -0.35 
\put{$($} at 11.7 0
\put{$)$} at 13.3 0 
\end{tiny}

\circulararc 45 degrees from 11.25 1.5  
center at 7.25 -8

\plot 
11.25 1.5 
11.03 1.5 / 

\plot 
11.25 1.5 
11.085 1.67 /

\put{Figure 2} at 7.12 -1.5

\put{} at -1 0

\endpicture

%%%%%%%%%%%%%%%%%%%%%%%%%%%%%%%%%%%%%%%%%%%%%%%%%%%%%%%%%%%%%%%%%%%%%%%%%%%%%%%%%%%%%%%%%%%%%%

\vspace{1cm}

Now for (ii) put $A = |[u_{k+1}^i,v_{k+1}^i]|$, 
$B = |[a_{2^{k+1}-i},a_{2^{k+1}-i-1}]|$, $C = |[u_{k+1}^{i+1},v_{k+1}^{i+1}]|$, 
and $D = |[a_{2^{k+1}-i-1},a_{2^{k+1}-i-2}]|$. Then 
$$\left| \frac{|[a_{2^{k+1} - i - 1},u_{k+1}^{i+1}]|}{|[a_{2^{k+1} - i}, u_{k+1}^i]|} 
- 1 \right| \frac{1}{|[a_{2^{k+1} - i},u_{k+1}^{i}]|^{\alpha}} = 
\left| \frac{D - C}{B - A} - 1 \right| \frac{2^{\alpha}}{(B - A)^{\alpha}}.$$
Moreover, since \esp $A \leq B/2$ \esp and \esp $C = \lambda_k A$, \esp 
\begin{eqnarray*}
\left| \frac{D - C}{B - A} - 1 \right| &\leq& \left| \frac{D - B}{B - A} \right|  
+ \left| \frac{C - A}{B - A} \right| \leq 2 \left| \frac{D - B}{B} \right|  
+ \left| \lambda_k - 1 \right| \\
&=& \frac{M}{B} \left[ \frac{1}{(2^{k+1}-i-2)^{1+\varepsilon}} - 
\frac{1}{(2^{k+1}-i-1)^{1+\varepsilon}} \right] + M \frac{k}{2^k}\\ 
&\leq& M B \left[ (2^{k+1}-i-1)^{1+\varepsilon} - (2^{k+1}-i-2)^{1+\varepsilon} \right] 
+ M \frac{k}{2^k}\\
&\leq& \frac{M}{2^{k (1 + \varepsilon)}} 2^{k \varepsilon} + M \frac{k}{2^k}\\ 
&\leq& M \frac{k}{2^k}.
\end{eqnarray*}
%A similar argument shows that
%$$\left| \frac{D}{B} - 1 \right| \leq M B^{1/(1 + \varepsilon)}.$$
Therefore,
$$\left| \frac{D - C}{B - A} - 1 \right| \frac{2^{\alpha}}{(B - A)^{\alpha}} 
\leq M \frac{k}{2^k} 2^{k (1 + \varepsilon) \alpha},$$ 
and thus 
\begin{equation}
\left| \frac{|[a_{2^{k+1} - i - 1},u_{k+1}^{i+1}]|}{|[a_{2^{k+1} - i}, u_{k+1}^i]|} 
- 1 \right| \frac{1}{|[a_{2^{k+1} - i},u_{k+1}^{i}]|^{\alpha}} 
\leq M \frac{k}{2^{k (1 - (1 + \varepsilon) \alpha)}}.
\label{segunda}
\end{equation}
Finally, notice that by construction the estimates for (iii) are the same as those for (ii).

\vspace{0.25cm}

%%%%%%%%%%%%%%%%%%%%%%%%%%%%%%%%%%%%%%%%%%%%%%%%%%

\noindent{\bf Estimates for $g$:} The diffeomorphism $g$ is obtained 
by fitting together many maps sending:

\vspace{0.1cm}

\noindent (i) $[b_k,u_k^0]$ \esp\esp into \esp \esp $[b_k,v_k^0]$,

\vspace{0.1cm}

\noindent (ii) $[u_k^0,c_k]$ \esp\esp into \esp\esp $[v_k^0,c_k]$,

\vspace{0.1cm}

\noindent (iii) $[a_{2^k + 1},b_k]$ and $[c_k,a_{2^k}]$ into themselves as the identity. 
%(thus, there is no need of an estimate here).

\vspace{0.1cm}

For (i) notice that
$$\left| \frac{|[b_k,v_k^0]|}{[b_k,u_k^0]} - 1 \right| \frac{1}{|[b_k,u_k^0]|^{\alpha}} 
= \frac{|[u_k^0,v_k^0]|}{|[b_k,u_k^0]|^{1 + \alpha}} \leq 
\frac{2^{1 + \alpha} |[u_k^0,v_k^0]|}{\left( |[b_k,c_k]| - |[u_k^0,v_k^0]| \right)^{1 + \alpha}} = 
\frac{2^{1 + \alpha} |[b_k,c_k]|^{1 + \theta}}{(|[b_k,c_k]| - |[b_k,c_k]|^{1 + \theta})^{1 + \alpha}},$$
and thus
\begin{equation}
\left| \frac{|[b_k,v_k^0]|}{[b_k,u_k^0]} - 1 \right| \frac{1}{|[b_k,u_k^0]|^{\alpha}} 
\leq M |[b_k,c_k]|^{\theta - \alpha}.
\label{tercera}
\end{equation}

Finally, the estimates for (ii) are similar to those for (i) and we leave them 
to the reader.

\vspace{0.25cm}

%%%%%%%%%%%%%%%%%%%%%%%%%%%%%%%%%%%%%%%%%%%%%%%%%%

\noindent{\bf The choice of the parameters: }  According to our Lemma, and due to 
(\ref{primera}), (\ref{segunda}), and (\ref{tercera}), sufficient conditions for 
the $\ce^{1 + \alpha}$ smoothness of $f,g$ are:

\vspace{0.1cm}

\noindent -- \esp \esp $(1 + \varepsilon) (1 + \theta) \alpha < 1$,

\vspace{0.1cm}

\noindent -- \esp \esp $\frac{1}{1 + \varepsilon} > \alpha$

\vspace{0.1cm}

\noindent -- \esp \esp $\theta > \alpha$.

\vspace{0.1cm}

Now, for \esp $0 < \alpha < (\sqrt{5} - 1) / 2$ \esp one easily checks that these 
conditions are satisfied for \esp $\theta = \alpha + \varepsilon$, \esp where 
\esp $\varepsilon > 0$ \esp is small enough so that \esp 
$(1 + \varepsilon) (1 + \alpha + \varepsilon) \alpha < 1$. 

\vspace{0.3cm}

\noindent{\bf Acknowledgments.} I would like to thank J. Cantwell and L. Conlon 
for motivating me to work on and write out the example of this Note.

%%%%%%%%%%%%%%%%%%%%%%%%%%%%%%%%%%%%%%%%%%%%%%%%%%%%%%%%%%%%%%%%%%%%%%%%%%%%%%%%%%%%%%%%%%%%%%%%

\begin{small}

\vspace{0.1cm}

\noindent Andr\'es Navas, Univ. de Santiago, Alameda 3363, Chile 
(anavas@usach.cl)\\

\end{small}

\end{document}